\documentclass{amsart}

\usepackage{amsthm}
\usepackage{amsmath}
\usepackage{amssymb}
\usepackage{amsfonts}
\usepackage{latexsym}
\usepackage[all]{xy} \SelectTips{eu}{}
\usepackage{hyperref}
\usepackage{eucal}
\usepackage{xcolor}
\newtheorem{intthm}{Theorem}[]

\newtheorem*{intque*}{Question}
\newtheorem*{intexa*}{Example}
\newcommand{\numberseries}{\bfseries}   %Fontseries used for numbering
\newlength{\thmtopspace}                %Space above theorem
\newlength{\thmbotspace}                %Space below theorem
\newlength{\thmheadspace}               %Space after theorem label
\newlength{\thmindent}                  %For indenting
\newtheoremstyle{bfupright head,slanted body}
{\thmtopspace}{\thmbotspace}
{\slshape}{\thmindent}{\bfseries}{.}{\thmheadspace}
{{\numberseries \thmnumber{#2\;}}\thmnote{#3}}

\newtheoremstyle{bfupright head,upright body}
{\thmtopspace}{\thmbotspace}
{\upshape}{\thmindent}{\bfseries}{.}{\thmheadspace}
{{\numberseries \thmnumber{#2\;}}\thmnote{#3}}

\newtheoremstyle{fixed bf head,slanted body}
{\thmtopspace}{\thmbotspace}{\slshape}
{\thmindent}{\bfseries}{.}{\thmheadspace}
{{\numberseries \thmnumber{#2\;}}\thmname{#1}\thmnote{ (#3)}}

\newtheoremstyle{fixed bf head,upright body}
{\thmtopspace}{\thmbotspace}{\upshape}
{\thmindent}{\bfseries}{.}{\thmheadspace}
{{\numberseries \thmnumber{#2\;}}\thmname{#1}\thmnote{ (#3)}}

\newtheoremstyle{numbered paragraph}
{\thmtopspace}{\thmbotspace}{\upshape}
{\thmindent}{\upshape}{}{\thmheadspace}
{{\numberseries \thmnumber{#2.}}}

\theoremstyle{bfupright head,slanted body}
\newtheorem{res}{}[section]             \newtheorem*{res*}{}

\theoremstyle{bfupright head,upright body}
\newtheorem{bfhpg}[res]{}               \newtheorem*{bfhpg*}{}

\theoremstyle{fixed bf head,slanted body}
\newtheorem{thm}[res]{Theorem}          \newtheorem*{thm*}{Theorem}
\newtheorem{prp}[res]{Proposition}      \newtheorem*{prp*}{Proposition}
\newtheorem{cor}[res]{Corollary}        \newtheorem*{cor*}{Corollary}
\newtheorem{lem}[res]{Lemma}            \newtheorem*{lem*}{Lemma}
         \newtheorem*{que*}{Question}
\theoremstyle{fixed bf head,upright body}
\newtheorem{dfn}[res]{Definition}       \newtheorem*{dfn*}{Definition}
\newtheorem{rmk}[res]{Remark}           \newtheorem*{rmk*}{Remark}
           \newtheorem*{exa*}{Example}
\newtheorem{setup}[res]{Setup}           \newtheorem*{setup*}{Setup}
\newtheorem{nota}[res]{Notation}           \newtheorem*{nota*}{Notation}
\theoremstyle{numbered paragraph}
\newtheorem{ipg}[res]{}
\newlength{\thmlistleft}        %leftmargin
\newlength{\thmlistright}       %rightmargin
\newlength{\thmlistpartopsep}   %partopsep
\newlength{\thmlisttopsep}      %topsep
\newlength{\thmlistparsep}      %parsep
\newlength{\thmlistitemsep}     %itemsep

\newcounter{eqc}
	{\end{list}}%

\newcounter{prt}
\newenvironment{prt}{\begin{list}{\upshape (\alph{prt})}%
		{\usecounter{prt}%
			\setlength{\leftmargin}{\thmlistleft}%
			\setlength{\labelwidth}{\thmlistleft}%
			\setlength{\rightmargin}{\thmlistright}%
			\setlength{\partopsep}{\thmlistpartopsep}%
			\setlength{\topsep}{\thmlisttopsep}%
			\setlength{\parsep}{\thmlistparsep}%
			\setlength{\itemsep}{\thmlistitemsep}}}%
	{\end{list}}%

\newcounter{rqm}
\newenvironment{rqm}{\begin{list}{\upshape (\arabic{rqm})}%
		{\usecounter{rqm}%
			\setlength{\leftmargin}{\thmlistleft}%
			\setlength{\labelwidth}{\thmlistleft}%
			\setlength{\rightmargin}{\thmlistright}%
			\setlength{\partopsep}{\thmlistpartopsep}%
			\setlength{\topsep}{\thmlisttopsep}%
			\setlength{\parsep}{\thmlistparsep}%
			\setlength{\itemsep}{\thmlistitemsep}}}%
	{\end{list}}%

	{\end{list}}%

\newenvironment{prf*}[1][Proof]{%
	\begin{proof}[\bf #1]
		\setcounter{equation}{0}
		}
	{\end{proof}
}

\newcommand{\pgref}[1]{\ref{#1}}

\renewcommand{\eqref}[1]{(\pgref{eq:#1})}

\newcommand{\thmcite}[2][?]{\cite[thm.~#1]{#2}}
\newcommand{\prpcite}[2][?]{\cite[prop.~#1]{#2}}

\newcommand{\corcite}[2][?]{\cite[cor.~#1]{#2}}
\newcommand{\lemcite}[2][?]{\cite[lem.~#1]{#2}}

\newcommand{\dfncite}[2][?]{\cite[def.~#1]{#2}}

\newcommand{\exacite}[2][?]{\cite[exa.~#1]{#2}}

\numberwithin{equation}{res}
\def\urltilda{\kern -.15em\lower .7ex\hbox{\~{}}\kern .04em}
\newcommand{\setof}[3][\mspace{1mu}]{\{#1#2 \mid #3#1\}}
\newcommand{\GF}[1]{\mathsf{GF}(#1)}
\newcommand{\GP}[1]{\mathsf{GP}(#1)}
\newcommand{\PGF}[1]{\mathsf{PGF}(#1)}
\newcommand{\GI}[1]{\mathsf{GI}(#1)}
\newcommand{\Flat}[1]{\mathsf{Flat}(#1)}
\newcommand{\Cot}[1]{\mathsf{Cot}(#1)}
\newcommand{\Mod}[1]{{_{#1}\mathsf{Mod}}}
\newcommand{\Prj}[1]{\mathsf{Prj}(#1)}
\newcommand{\Inj}[1]{\mathsf{Inj}(#1)}
\newcommand{\ch}[1]{\mathsf{Ch}(#1)}
\newcommand{\dgPrj}[1]{\mathsf{dgPrj}(#1)}
\newcommand{\dgInj}[1]{\mathsf{dgInj}(#1)}
\newcommand{\dgFlat}[1]{\mathsf{dgFlat}(#1)}
\newcommand{\dgCot}[1]{\mathsf{dgCot}(#1)}
\newcommand{\D}[1]{\mathsf{D}(#1)}

\newcommand{\Hom}[3][]{\operatorname{Hom}_{#1}(#2,#3)}

\newcommand{\Ext}[4][]{\operatorname{Ext}_{#1}^{#2}(#3,#4)}
\newcommand{\Tor}[4][]{\operatorname{Tor}_{#2}^{#1}(#3,#4)}
\newcommand{\Coker}[1]{\nobreak{\operatorname{Coker}#1}}
\newcommand{\A}{\mathsf{A}}
\newcommand{\Acy}{\mathsf{Acy}}
\newcommand{\CC}{\mathsf{C}}
\newcommand{\E}{\mathsf{E}}
\newcommand{\FF}{\mathsf{F}}
\newcommand{\G}{\mathsf{G}}
\newcommand{\Se}{\mathsf{S}}
\newcommand{\T}{\mathsf{T}}
\newcommand{\R}{\mathsf{R}}
\newcommand{\W}{\mathsf{W}}

\newcommand{\colim}{\mbox{\rm colim}}
\newcommand{\C}{\mathcal{C}}
\newcommand{\Q}{\mathcal{Q}}
\newcommand{\pd}{\mbox{\rm pd}}
\newcommand{\fd}{\mbox{\rm fd}}
\newcommand{\id}{\mathrm{id}}
\newcommand{\kk}{\Bbbk}
\newcommand{\is}{\cong}

\newcommand{\lMod}[1]{{}_{#1}\mspace{-1mu}\operatorname{Mod}}
\newcommand{\rMod}[1]{\mathsf{Mod}_{#1}}

\def\soft#1{\leavevmode\setbox0=\hbox{h}\dimen7=\ht0\advance
	\dimen7 by-1ex\relax\if t#1\relax\rlap{\raise.6\dimen7
		\hbox{\kern.3ex\char'47}}#1\relax\else\if T#1\relax
	\rlap{\raise.5\dimen7\hbox{\kern1.3ex\char'47}}#1\relax
	\else\if d#1\relax\rlap{\raise.5\dimen7\hbox{\kern.9ex
			\char'47}}#1\relax\else\if D#1\relax\rlap{\raise.5\dimen7
		\hbox{\kern1.4ex\char'47}}#1\relax\else\if l#1\relax
	\rlap{\raise.5\dimen7\hbox{\kern.4ex\char'47}}#1\relax
	\else\if L#1\relax\rlap{\raise.5\dimen7\hbox{\kern.7ex
			\char'47}}#1\relax\else\message{accent \string\soft
		\space #1 not defined!}#1\relax\fi\fi\fi\fi\fi\fi}

\begin{document}

\title[Flat model structures and Gorenstein objects in functor categories]%
{Flat model structures and Gorenstein objects in functor categories}

\author[Z.X. Di]{Zhenxing Di}

\address{Zhenxing Di: School of Mathematical Sciences, Huaqiao University, Quanzhou 362021, China}

\email{dizhenxing@163.com}

\author[L.P. Li]{Liping Li}

\address{Liping Li: Department of Mathematics, Hunan Normal University, Changsha 410081, China}

\email{lipingli@hunnu.edu.cn}

\author[L. Liang]{Li Liang}

\address{Li Liang: (1) Department of Mathematics, Lanzhou Jiaotong University, Lanzhou 730070, China; (2) Gansu Provincial Research Center for Basic Disciplines of
Mathematics and Statistics, Lanzhou 730070, China}

\email{lliangnju@gmail.com}

\urladdr{https://sites.google.com/site/lliangnju}

\author[Y.J. Ma]{Yajun Ma}

\address{Yajun Ma: Department of Mathematics, Lanzhou Jiaotong University, Lanzhou 730070, China}

\email{13919042158@163.com}

\thanks{Z.X. Di was partly supported by NSF of China (Grant No. 11971388), the Scientific Research Funds of Huaqiao University (Grant No. 605-50Y22050) and the Fujian Alliance of Mathematics (Grant No. 2024SXLMMS04); L.P. Li was partly supported by NSF of China (Grant No. 12171146); L. Liang was partly supported by NSF of China (Grant No. 12271230) and the Foundation for Innovative Fundamental Research Group Project of Gansu Province (Grant No. 23JRRA684); Y. Ma was partly supported by NSF of Gansu Province (Grant No. 23JRRA866) and the Youth Foundation of Lanzhou Jiaotong University (Grant No. 2023023).}

%\date{\today}



\begin{abstract}
We construct a flat model structure on the category $\Mod{\Q,R}$ of additive functors from a small preadditive category $\Q$ satisfying certain conditions to the module category $\Mod{R}$ over an associative ring $R$, whose homotopy category is the $\Q$-shaped derived category introduced by Holm and J{\o}rgensen. Moreover, we prove that for an arbitrary associative ring $R$, an object in $\Mod{\Q,R}$ is Gorenstein projective (resp., Gorenstein injective, Gorenstein flat, projective coresolving Gorenstein flat) if and only if so is its value on each object of $\Q$, and hence improve a result by Dell'Ambrogio, Stevenson and \v{S}\v{t}ov\'{\i}\v{c}ek.
\end{abstract}

\maketitle

\thispagestyle{empty}
\section*{Introduction}%%%%%%%%%%%%%%%%%%%%%%%%%%%%%%%%%%
\noindent
Model category theory was introduced by Quillen in \cite{Qui67} to provide a uniform way of formally introducing homotopy theory into general categories, and was heavily applied in various areas such as algebraic topology, category theory, and representation theory; see for instance Hovey \cite{Ho99}. For representation theory, abelian model structures are of significant interest: people try to construct various abelian model structures on abelian categories to obtain the corresponding homotopy theories. The construction of abelian model structures turns out to be equivalent to the construction of Hovey triples by Hovey's correspondence established in \cite{Ho02}: an abelian model structure on an abelian category $\A$ is equivalent to a triple $(\CC, \W, \FF)$ of subcategories of $\A$ such that $\W$ is thick, and $(\CC, \W\cap\FF)$ and $(\CC\cap\W,\FF)$ are two complete cotorsion pairs, where $\W$ (resp., $\CC$ and $\FF$) is the subcategory of $\A$ consisting of all trivial (resp., cofibrant and fibrant) objects associated to the corresponding abelian model structure.

For an arbitrary associative ring $R$, \v{S}aroch and \v{S}t'ov\'{\i}\v{c}ek described in \cite{SS20} three non-trivial abelian model structures on $\Mod{R}$, the category of left $R$-modules:
\begin{itemize}
\item a projective model structure $(\PGF{R}, \PGF{R}^\perp, \Mod{R})$,
\item an injective model structure $(\Mod{R}, {^\perp\GI{R}}, \GI{R})$,
\item and a flat model structure $(\GF{R}, \PGF{R}^\perp, \Cot{R})$.
\end{itemize}
Here we recall that an abelian model structure $(\CC, \W, \FF)$ is called \emph{projective} (resp., \emph{flat}) if the subcategory $\CC\cap\W$ consists of projectives (resp., flats), and it is called \emph{injective} if the subcategory $\W\cap\FF$ consists of injectives.

When one considers the category $\ch{R}$ of chain complexes of left $R$-modules, there are more abelian model structures. In particular, there are
\begin{itemize}
\item a projective model structure $(\dgPrj{R}, \Acy, \ch{R})$,
\item an injective model structure $(\ch{R}, \Acy, \dgInj{R})$,
\item and a flat model structure $(\dgFlat{R}, \Acy, \dgCot{R})$,
\end{itemize}
and the corresponding homotopy categories of the above three model structures all coincide with $\D{R}$, the derived category of $R$; see \cite[exas. 3.2 and 3.3]{Ho02} and Gillespie \corcite[5.1]{Gi04}. Here $\Acy$ is the subcategory of acyclic complexes, and $\dgPrj{R}$ (resp., $\dgInj{R}$, $\dgFlat{R}$ and $\dgCot{R}$) is the subcategory of DG-projective (resp., DG-injective, DG-flat and DG-cotorsion) complexes.

Since $\ch{R}$ can be viewed as a functor category in a natural way (see \exacite[8.13]{HJ21}), one may wonder to extend the above model structures to functor categories. This was considered by Holm and J{\o}rgensen recently in \cite{HJ21}, where they constructed the projective and injective model structures on the category $\Mod{\Q,R}$ of additive functors from a small preadditive category $\Q$ satisfying certain conditions specified in \cite[Setup 2.5]{HJ21} (see also Setup \ref{setup}) to the category $\Mod{R}$ for any associative ring $R$. Explicitly, they proved in \thmcite[A]{HJ21} that there are
\begin{itemize}
\item a projective model structure $({^\perp\E}, \E, \Mod{\Q,R})$,
\item and an injective model structure $(\Mod{\Q,R}, \E, \E^\perp)$
\end{itemize}
on $\Mod{\Q,R}$, where $\E = \{X \in \Mod{\Q,R}\ |\ X^{\natural}\ \mathrm{has\ finite\ projective\ dimension\ in}\ \Mod{\Q}\}$. Here $(-)^{\natural}$ denotes the forgetful functor from $\Mod{\Q,R}$ to $\Mod{\Q}$. Both model structures have the same homotopy category called the $\Q$-shaped derived category. In particular, by taking $\Q$ to be a special category, their result includes the projective and injective model structures on $\ch{R}$ as examples.

We note that the flat model structure on $\ch{R}$ was not extended to the more general framework in \cite{HJ21}. As the first main result of this paper, we fill this gap by describing a flat model structure on $\Mod{\Q,R}$ whose homotopy category also coincides with the $\Q$-shaped derived category. The following is a special case of Corollary \ref{flat model structure} with $\kk=\mathbb{Z}$.

\begin{intthm}\label{thmA}
Let $\Q$ be a small preadditive category satisfying the conditions specified in Setup \ref{setup}. Then for any associative ring $R$, there is a hereditary abelian model structure
\[
(^{\perp}(\Cot{\Q,R}\cap\E), \E, \Cot{\Q,R})
\]
on $\Mod{\Q,R}$, where $^{\perp}(\Cot{\Q,R}\cap\E)$ is the subcategory of cofibrant objects, $\E$ is the subcategory of trivial objects, and $\Cot{\Q,R}$ is the subcategory of fibrant objects.
\end{intthm}

Members in $\Cot{\Q,R}$ are cotorsion objects in $\Mod{\Q,R}$ defined in \ref{flat-cotorsion}.

%\begin{rmk*}
%Actually, in the paper we not only provide a flat model structure on $\Mod{\Q,A}$, but also generalize Holm and J{\o}rgensen's projective/injective model structures from the following two aspects: firstly, we only assume that $\kk$ is a commutative ring and drop the condition that $\kk$ is Gorenstein imposed in \cite{HJ21}; secondly, $\Q$ only needs to satisfy the condition that the $\kk$-module $\Q(p,q)$ is projective for each pair of objects $p$ and $q$. For details, the reader may compare Propositions \ref{prj model} and \ref{inj model} to the corresponding results in \cite{HJ21}.
%\end{rmk*}

%\begin{equation*}
%  \ast \ \ \ast \ \ \ast
%\end{equation*}
We mention that the key ingredient to prove Theorem \ref{thmA} is a recent result on Gorenstein homological objects by \v{S}aroch and \v{S}t'ov\'{\i}\v{c}ek \cite{SS20}; see Lemma \ref{known}. For a Gorenstein ring $R$ and a small preadditive category $\Q$ satisfying the conditions (1)-(3) in Setup \ref{setup}, it is known that an object $X$ in $\Mod{\Q, R}$ is Gorenstein projective if and only if each left $R$-module $X(q)$ is Gorenstein projective for $q\in\Q$; see Dell'Ambrogio, Stevenson and \v{S}\v{t}ov\'{\i}\v{c}ek \corcite[4.8]{DSS17}. Note that chain complexes of left $R$-modules can be regarded as additive functors from a special small preadditive category to the category $\Mod{R}$ of left $R$-modules. The result \corcite[4.8]{DSS17} is a generalization of \thmcite[4.5]{EEnJGR98} by Enochs and Garc{\'\i}a Rozas: an unbounded chain complex $M$ of left $R$-modules is Gorenstein projective if and only if each left $R$-module $M_i$ is Gorenstein projective for $i \in \mathbb{Z}$ under the extra assumption that $R$ is a Gorenstein ring. Since in \cite{LZ09,YL11} Liu, Yang and Zhang have shown that \thmcite[4.5]{EEnJGR98} actually holds for arbitrary associative rings, it is then natural to ask if \corcite[4.8]{DSS17} also holds for these small preadditive categories $\Q$ without any restriction on the coefficient ring $R$.

Our second main result answers this question affirmatively. Relying on the concept of Frobenius functors introduced by Chen and Ren in \cite{CR22} and a main result of that paper, we prove that an additive functor $X$ from $\Q$ to $\Mod{R}$ is Gorenstein projective (resp., Gorenstein injective, Gorenstein flat, projective coresolving Gorenstein flat) if and only if so is each left $R$-module $X(q)$ in the following result, which is a special case of Theorems \ref{GP}, \ref{GI}, \ref{GF} and \ref{PGF} with $\kk=\mathbb{Z}$.

\begin{intthm}\label{thmB}
Let $\Q$ be a small preadditive category satisfying the conditions (1)-(3) in Setup \ref{setup} and $R$ an arbitrary associative ring. Then an object $X\in \Mod{\Q,R}$ is Gorenstein projecitve (resp., Gorenstein injective, Gorenstein flat, projectively coresolved Gorenstein flat) if and only if each left $R$-module $X(q)$ is Gorenstein projecitve (resp., Gorenstein injective, Gorenstein flat, projectively coresolved Gorenstein flat) for $q\in\Q$.
\end{intthm}

\section{Preliminaries}\label{property}
\noindent
Throughout this paper, by the term \emph{``subcategory''} we always mean a full subcategory which is closed under isomorphisms and contains the zero object. We assume that $\kk$ is a commutative ring and $A$ is an associative $\kk$-algebra, and assume that $\C$ is a small $\kk$-linear category such that $\C(p,q)$ is a projective $\kk$-module for each pair of objects $p$ and $q$ in $\C$.

Let $\Mod{\C,A}$ be the category of $\kk$-linear functors from $\C$ to the left $A$-module category $\Mod{A}$, and let $\rMod{\C,A}$ be the category of $\kk$-linear functors from $\C^{\sf op}$ to the right $A$-module category $\rMod{A}$. In this section, we collect some basic results on homological objects in $\Mod{\C,A}$, which will be used frequently in the paper.

\begin{bfhpg}[\bf Cotorsion pairs]\label{cotpair}
Let $\A$ be an abelian category with enough projectives and injectives.
A pair $(\CC, \W)$ of subcategories of $\A$ is called a \emph{cotorsion pair} if $\CC^{\bot}=\W$ and $^{\bot}\W=\CC$. Here
$$\quad \ \ {\CC}^{\bot}=\{N\in {\sf A}\ |\ {\rm Ext}_{\A}^{1}(C,N)=0 {\rm \ for\ all}\ C\in\CC\},\ \text{and}$$
$$^{\bot}{\W}=\{M\in {\sf A}\ |\ {\rm Ext}_{\A}^{1}(M,W)=0 {\rm \ for\ all}\ W\in\W\}.$$
It is clear that $({^\perp(\CC^\perp)},\CC^\perp)$ is a cotorsion pair, which is called a cotorsion pair \emph{generated} by $\CC$. Recall from e.g. Enochs and Jenda \cite{rha} that a cotorsion pair $(\CC,\W)$ is \emph{complete} if for each object $M$ in $\A$ there is an exact sequence $0 \to M \to W \to C \to 0$ with $W\in\W$ and $C\in\CC$ (or there is an exact sequence $0 \to W' \to C' \to M \to 0$ with $W'\in\W$ and $C'\in\CC$). A cotorsion pair $(\CC,\W)$ is called \emph{hereditary} if $\CC$ is closed under kernels of epimorphisms or $\W$ is closed under cokernels of monomorphisms.
\end{bfhpg}

The following notion of projective/injective cotorsion pairs can be found in Gillespie \dfncite[3.4]{Gil16}.

\begin{dfn}
Let $\A$ be an abelian category with enough projectives and injectives. A complete cotorsion pair $(\CC, \W)$ is said to be \emph{projective} if $\W$ is thick (that is, if two out of three of the terms in a short exact sequence are in $\W$, then so is the third) and $\CC\cap\W$ coincides with the subcategory of projective objects. The notion of \emph{injective} cotorsion pairs is defined dually.
\end{dfn}

\begin{rmk}
A projective cotorsion pair $(\CC, \W)$ determines a projective model structure on $\A$ such that $\CC$ is the subcategory of cofibrant objects, $\W$ is the subcategory of trivial objects, and all objects in $\A$ are fibrant objects. Dually, an injective cotorsion pair $(\W, \FF)$ determines an injective model structure on $\A$ such that $\FF$ is the subcategory of fibrant objects, $\W$ is the subcategory of trivial objects, and all objects in $\A$ are cofibrant objects.
\end{rmk}

\begin{bfhpg}[\bf Flat objects]\label{representable}
Let $\G$ be a locally finitely presentable (in the sense of Crawley-Boevey \cite[\S1]{WCB94}) Grothendieck category admitting enough projectives. A sequence
\[
0 \to K\to M\to N\to 0
\]
in $\G$ is said to be \emph{pure} if for each finitely presentable object $T$, the sequence
\[
0\to \Hom[\G]{T}{K} \to \Hom[\G]{T}{M} \to \Hom[\G]{T}{N} \to 0
\]
is exact (in this case, $M\to N\to 0$ is called a pure epimorphism). An object $F$ in $\G$ is called \emph{flat} if every epimorphism $M\to F$ is a pure epimorphism. It is not difficult to check that $F$ is flat if and only if it is a filtered colimit of finitely generated projective objects. Let ${\sf Flat(\G)}$ be the subcategory of flat objects in $\G$. Then ${\sf Flat(\G)}$ is closed under extensions, filtered colimits and pure epimorphisms by Rump \prpcite[7]{Ru10}, and so the pair $({\sf Flat(\G)}, {\sf Flat(\G)}^\perp)$ is a complete and hereditary cotorsion pair generated by a set by the proof of \prpcite[A.3]{HJ21}, that is, there is a set $\Se\subseteq{\sf Flat(\G)}$ such that $\Se^\perp = {\sf Flat(\G)}^\perp$.
\end{bfhpg}

\begin{ipg}\label{flat-cotorsion}
It follows from \prpcite[3.12]{HJ21} that $\Mod{\C,A}$ is a locally finitely presentable Grothendieck category admitting enough projectives. Let $\Flat{\C,A}$ be the subcategory of flat objects in $\Mod{\C,A}$. An object $X$ in $\Mod{\C,A}$ is called \textit{cotorsion} if it is in $\Flat{\C,A}^\perp$. Denote the subcategory of cotorsion objects by $\Cot{\C,A}$. Then by \ref{representable}, the pair $(\Flat{\C,A}, \Cot{\C,A})$ is a complete and hereditary cotorsion pair generated by a set.
\end{ipg}

\begin{bfhpg}[\bf Tensor products in $\Mod{\C,A}$]\label{tensor}
Following Oberst and R\"{o}hrl \cite{Oberst70}, we define a tensor product functor $-\otimes_{\C,A}-: \rMod{\C,A}\times \Mod{\C,A}\to \Mod{\kk}$; it is right exact in each variable.

Let $X\in \Mod{\C,A}$. We construct for any $\Bbbk$-module $G$, an object $\Hom[\Bbbk]{X}{G}\in \rMod{\C,A}$ as follows:
 \begin{itemize}
\item For each $q\in\C^{\sf op}$ let $\Hom[\Bbbk]{X}{G}(q)=\Hom[\Bbbk]{X(q)}{G}\in \rMod{A}$.
\item For each morphism $f: p\to q$ in $\C^{\sf op}$ let $\Hom[\Bbbk]{X}{G}(f)=\Hom[\Bbbk]{X(f^{\sf op})}{G}:\Hom[\Bbbk]{X(p)}{G}\to \Hom[\Bbbk]{X(q)}{G}$.
\end{itemize}
It is evident that $\Hom[\Bbbk]{X}{-}$ is a functor from the category $\Mod{\Bbbk}$ to $\rMod{\C,A}$. This functor is left exact and preserves arbitrary products, so it has a left adjoint functor from $\rMod{\C,A}$ to $\Mod{\Bbbk}$,
which is denoted $-\otimes_{\C,A} X$ and will play the role of the tensor product in our work. Then for each $Y\in\rMod{\C,A}$ there are natural isomorphisms
\begin{equation*}\label{1.6.1}
\tag{\ref{tensor}.1}
\Hom[\Bbbk]{Y\otimes_{\C,A} X}{G}\cong \Hom[\C^{\sf op},A^{\sf op}]{Y}{\Hom[\Bbbk]{X}{G}}\cong\Hom[\C,A]{X}{\Hom[\kk]{Y}{G}}.
\end{equation*}
\end{bfhpg}

\begin{nota}
Let $\Prj{\C,A}$ (resp., $\Inj{\C,A}$) be the subcategory of projective objects (resp., injective objects) in $\Mod{\C,A}$. In the case that $A=\kk$, we write $\Mod{\C}$ (resp., $\Prj{\C}$, $\Inj{\C}$, $\Flat{\C}$ and $\Cot{\C}$) instead of $\Mod{\C,\kk}$ (resp., $\Prj{\C,\kk}$, $\Inj{\C,\kk}$, $\Flat{\C,\kk}$ and $\Cot{\C,\kk}$), and we write $-\otimes_{\C}-$ instead of $-\otimes_{\C,\kk}-$.
\end{nota}

\begin{ipg}\label{adjoint triple}
By \cite[Corollary 3.9]{HJ21}, for every object $q$ in $\C$ there is an adjoint triple $(F_{q},E_{q},G_{q})$ as follows:
  \begin{equation*}
  \xymatrix@C=4pc{
    \Mod{\C,A}
    \ar[r]^-{E_{q}}
    &
    \Mod{A}
    \ar@/_1.8pc/[l]_-{F_{q}}
    \ar@/^1.8pc/[l]^-{G_{q}}
  }
  \qquad \text{given by} \qquad
  {\setlength\arraycolsep{1.5pt}
   \renewcommand{\arraystretch}{1.2}
  \begin{array}{rcl}
  F_q(M) &=& \C(q,-) \otimes_\kk M \\
  E_q(X) &=& X\mspace{1.5mu}(q) \\
  G_{q}(M) &=& \Hom[\kk]{\C(-, q)}{M}\;.
  \end{array}
  }
  \end{equation*}
\end{ipg}

It is clear that the \textit{evaluation} functor $E_q$ is exact. Since $\C(p,q)$ is a projective $\kk$-module for each pair of objects $p$ and $q$ in $\C$ by the assumption, both the functors $F_q$ and $G_q$ are exact as well. Then by Holm and J{\o}rgensen \lemcite[5.1]{HJ19}, we have the following result.

\begin{lem}\label{n-iso}
Let $X$ be an object in $\Mod{\C,A}$ and $q$ an object in $\C$. Then for each $M \in \Mod{A}$ and $n \geqslant 0$, there are isomorphisms
$$\Ext[\C,A]{n}{F_q(M)}{X}\is\Ext[A]{n}{M}{X(q)}$$
and
$$\Ext[A]{n}{X(q)}{M}\is\Ext[\C,A]{n}{X}{G_q(M)}.$$
\end{lem}

The following proposition asserts that the evaluation functor preserves special homological property: if an object (which is a functor) in $\Mod{\C,A}$ has a certain special homological property, so does its value on each object in $\C$.

\begin{prp}\label{degreewise}
Let $X$ be an object in $\Mod{\C,A}$ and $q$ an object in $\C$. Then the following statements hold.
\begin{prt}
\item If $X\in\Prj{\C,A}$, then $X(q)\in\Prj{A}$.
\item If $X\in\Inj{\C,A}$, then $X(q)\in\Inj{A}$.
\item If $X\in\Flat{\C,A}$, then $X(q)\in\Flat{A}$.
\item If $X\in\Cot{\C,A}$, then $X(q)\in\Cot{A}$.
\end{prt}
\end{prp}

\begin{prf*}
Statements (a) and (b) hold since the functors $F_q$ and $G_q$ are exact. For $X \in \Flat{\C,A}$, one has $X \is \colim_iP_i$ with $P_i \in \Prj{\C,A}$, and deduce by (a) that
\[
X(q) = E_q(X) \is E_q(\colim_i P_i) \is \colim_i E_q(P_i) = \colim_i P_i(q) \in\Flat{A},
\]
so the statement (c) is true.

To verify the statement (d), we fix a flat left $A$-module $M$ and write $M = \colim_i T_i$ with each $T_i$ a projective left $A$-module. It follows that each $F_q(T_i)$ is projective and $F_q(M)\is\colim_iF_q(T_i)$ is flat in $\Mod{\C,A}$. By Lemma \ref{n-iso}, one has
\[
\Ext[A]{1}{M}{X(q)}\is \Ext[\C,A]{1}{F_q(M)}{X}=0.
\]
Therefore, $X(q)$ is a cotorsion left $A$-module.
\end{prf*}

\begin{cor}\label{cot}
Let $N$ be a cotorsion left $A$-module and $q$ an object in $\C$. Then $G_{q}(N)$ is a cotorsion object in $\Mod{\C,A}$.
\end{cor}

\begin{prf*}
Take an arbitrary flat object $X$ in $\Mod{\C,A}$. Then by Lemma \ref{n-iso}
\[
\Ext[\C,A]{1}{X}{G_q(N)} \is \Ext[A]{1}{X(q)}{N}=0,
\]
as $X(q)$ is flat by Proposition \ref{degreewise}(c).
\end{prf*}

\begin{ipg}
Let $(-)^{\natural}$ denote the forgetful functor from $\Mod{\C,A}$ to $\Mod{\C}$. Then $(-)^{\natural}$ is an exact functor and there is an adjoint triple as follows (see \corcite[3.5]{HJ21}):
  \begin{equation*}
  \xymatrix@C=4pc{
    \Mod{\C,A}
    \ar[r]^-{(-)^{\natural}}
    &
    \Mod{\C}
    \ar@/_1.8pc/[l]_-{-\otimes_{\kk}A}
    \ar@/^1.8pc/[l]^-{\Hom[\kk]{A}{-}}
  }
  \end{equation*}
\end{ipg}

\begin{prp}\label{forgive}
Let $X$ be an object in $\Mod{\C,A}$. Then the following statements hold.
\begin{prt}
\item If $X \in \Prj{\C,A}$, then the projective dimension of $X^\natural$ is bounded by $\pd_\kk A$.
\item If $X \in \Inj{\C,A}$, then the injective dimension $X^\natural$ is bounded by $\fd_\kk A$.
\item If $X \in \Flat{\C,A}$, then the flat dimension of $X^\natural$ is bounded by $\pd_\kk A$.
\end{prt}
\end{prp}
\begin{prf*}
In the following, we use the symbols $\pd_\C(-)$, $\id_\C(-)$ and $\fd_\C(-)$ to denote the projective, injective and flat dimension of objects in $\Mod{\C}$, respectively.

(a). For each object $q$ in $\C$, one has $F_q(A)^\natural = \C(q,-) \otimes_\kk A^\natural$. Since the functor
\[
\C(q,-) \otimes_\kk-: \Mod{\kk} \to \Mod{\C}
\]
preserves projective objects and is exact, we deduce that
\[
\pd_\C(F_q(A)^\natural) \leqslant \pd_\kk A^\natural = \pd_\kk A.
\]
Thus by \prpcite[3.12(a)]{HJ21}, one has $\pd_\C X^\natural \leqslant \pd_\kk A$.

(b). For each object $q$ in $\C$ and $I \in \Inj{A}$, one has $G_q(I)^\natural = \Hom[\kk]{\C(-,q)}{I^\natural}$. Since the functor $\Hom[\kk]{\C(-,q)}{-}: \Mod{\kk}\to \Mod{\C}$ preserves injective objects and is exact, one has
\[
\id_\C(G_q(I)^\natural) \leqslant \id_\kk I^\natural=\id_\kk I \leqslant \fd_\kk A,
\]
where the last inequality holds by Avramov and Foxby \corcite[4.2(b)(I)]{LLAHBF91}. It then follows from \prpcite[3.12(b)]{HJ21} that
$\id_\C X^\natural \leqslant \fd_\kk A$.

(c). Write $X = \colim_iP_i$ with each $P_i \in \Prj{\C,A}$. Then $X^\natural = \colim_iP_i^\natural$. By the statement (a), one has that $\pd_\C(P_i^\natural) \leqslant \pd_\kk A$ for each $i$, so
$\fd_\C(X^\natural) \leqslant \pd_\kk A$.
\end{prf*}

\begin{lem}\label{iso}
Assume that $A$ has finite projective dimension over $\kk$. Let $X$ be an object in $\Mod{\C,A}$. Then the following statements hold.
\begin{prt}
\item For each $G \in \Mod{\C}$ admitting an exact sequence $0 \to G \to F_{-1} \to F_{-2} \to \cdots$ with each $F_i \in \Flat{\C}$, there is an isomorphism
\[
\Ext[\C,A]{1}{G \otimes_\kk A}{X} \is \Ext[\C]{1}{G}{X^{\natural}}.
\]
\item For each $G \in \Mod{\C}$ admitting an exact sequence $\cdots \to I_1 \to I_0 \to G \to 0$ with each $I_i \in \Inj{\C}$, there is an isomorphism
\[
\Ext[\C,A]{1}{X}{\Hom[\kk]{A}{G}} \is \Ext[\C]{1}{X^{\natural}}{G}.
\]
\end{prt}
\end{lem}

\begin{prf*}
By Holm and J{\o}rgensen \lemcite[1.2]{HJ19b}, there is an exact sequence
\begin{equation*}
0 \to \Ext[\C,A]{1}{G \otimes_\kk A}{X} \to \Ext[\C]{1}{G}{X^{\natural}} \to \Hom[\C,A]{L_1(- \otimes_\kk A)(G)}{X}.
\end{equation*}
Statement (a) holds after we show that $L_1(-\otimes_\kk A)(G) = 0$. For this purpose, consider the exact sequence
\[
0 \to K \to P \to G \to 0
\]
in $\Mod{\C}$ with $P \in \Prj{\C}$. Clearly, $L_1(- \otimes_\kk A)(P) = 0$, so it is enough to show the exactness of the sequence
\begin{equation}\label{seq1}\tag{\ref{iso}.1}
0 \to K\otimes_\kk A \to P\otimes_\kk A \to G\otimes_\kk A \to 0
\end{equation}
in $\Mod{\C,A}$.

Consider the exact sequence
\[
0 \to G \to F_{-1} \to F_{-2} \to \cdots
\]
with each $F_i \in \Flat{\C}$. For each object $q$ in $\C$, the sequence
\[
0 \to G(q) \to F_{-1}(q) \to F_{-2}(q) \to \cdots
\]
is exact and each $F_i(q)$ is a flat $\kk$-module by Proposition \ref{degreewise}(c). Since $\fd_\kk A<\infty$ by the assumption, one gets that $\Tor[\kk]{1}{G(q)}{A} = 0$ by dimension shifting. Thus the sequence
\[
0 \to K(q) \otimes_\kk A \to P(q) \otimes_\kk A \to G(q) \otimes_\kk A \to 0
\]
is exact, which guarantees the exactness of the sequence (\ref{seq1}).

Statement (b) can be proved dually using \lemcite[1.3]{HJ19b} and Proposition \ref{degreewise}(b), observing that $\pd_\kk A<\infty$, and $((-)^{\natural}, \Hom[\kk]{A}{-})$ is an adjoint pair.
\end{prf*}

\begin{prp}\label{forgive-cotorsion}
Assume that $A$ has finite projective dimension over $\kk$. Let $X$ be an object in $\Mod{\C,A}$. If $X \in \Cot{\C,A}$, then $X^\natural \in \Cot{\C}$.
\end{prp}
\begin{prf*}
Take an object $F$ in $\Flat{\C}$ and write $F = \colim_iP_i$ with each $P_i \in \Prj{\C}$. Since the functor $- \otimes_{\kk} A$ preserves projectives, each $P_i \otimes_{\kk} A$ is contained in $\Prj{\C,A}$. Thus $F \otimes_{\kk} A$ belongs to $\Flat{\C,A}$, and by Lemma \ref{iso}(a) one deduces that
\[
\Ext[\C]{1}{F}{X^{\natural}} \is \Ext[\C,A]{1}{F \otimes_\kk A}{X} = 0
\]
because $X \in \Cot{\C,A}$. Consequently, $X^\natural$ is contained in $\Cot{\C}$.
\end{prf*}

\begin{bfhpg}[\bf Gorenstein homological objects]
Following \cite{rha}, an object $X$ in $\Mod{\C, A}$ is called \emph{Gorenstein projective} if there is an exact sequence
$$\cdots \to P_{1}\to P_{0}\to P_{-1} \to \cdots$$
in $\Mod{\C,A}$ with each $P_{i}$ projective such that $X\cong\Coker{(P_{1}\to P_{0})}$ and the sequence remains exact after applying the functor $\Hom[\C,A]{-}{P}$ for every projective object $P$ in $\Mod{\C,A}$.

Dually, an object $X$ in $\Mod{\C, A}$ is called \emph{Gorenstein injective} if there is an exact sequence
$$\cdots \to I_{1}\to I_{0}\to I_{-1}\to\cdots$$
in $\Mod{\C,A}$ with each $I_{i}$ injective such that $X\cong\Coker{(I_{1}\to I_{0})}$ and the sequence remains exact after applying the functor $\Hom[\C,A]{I}{-}$ for every injective object $I$ in $\Mod{\C,A}$.

An object $X$ in $\Mod{\C, A}$ is called \emph{Gorenstein flat} if there is an exact sequence
$$\cdots \to F_{1}\to F_{0}\to F_{-1}\to\cdots$$
in $\Mod{\C, A}$ with each $F_{i}$ flat such that $X\cong\Coker{(F_{1}\to F_{0})}$ and the sequence remains exact after applying the functor $E\otimes_{\C,A}-$ for every injective object $E$ in $\rMod{\C,A}$. Similarly, one can define \emph{projectively coresolved Gorenstein flat} objects in $\Mod{\C,A}$ when flat objects in the above exact sequence are replaced by projective objects; see \cite{SS20}.

Let $\GP{\C,A}$ (resp., $\GI{\C,A}$, $\GF{\C,A}$, $\PGF{\C,A}$) be the subcategory of Gorenstein projective objects (resp., Gorenstein injective objects, Gorenstein flat objects, projectively coresolved Gorenstein flat objects) in $\Mod{\C,A}$. In the case that $A=\kk$, we write $\GP{\C}$ (resp., $\GI{\C}$, $\GF{\C}$, $\PGF{\C}$) instead of $\GP{\C,A}$ (resp., $\GI{\C,A}$, $\GF{\C,A}$, $\PGF{\C,A}$).
\end{bfhpg}

\section{Model structures on $\Mod{\C,A}$}\label{model}
\noindent
In this section we describe projective, injective and flat model structures on the category $\Mod{\C,A}$, and in particular prove Theorem \ref{thmA} in the introduction. Throughout this section, we assume that $A$ is an associative $\kk$-algebra with $\pd_\kk A<\infty$.

The next result is from (the proofs of) \cite[thms. 4.9, 4.11 and 5.6, and cor 4.12]{SS20}; \v{S}aroch and \v{S}t'ov\'{\i}\v{c}ek proved it for modules over rings, but their techniques also work well for objects in $\Mod{\C}$ (see \cite[pp. 4]{SS20}).

\begin{lem}\label{known}
The following statements hold.
\begin{prt}
\item The pair $(\PGF{\C}, {\PGF{\C}}^\perp)$ is a complete hereditary cotorsion pair generated by a set such that ${\PGF{\C}}^\perp$ is thick and $\PGF{\C} \cap {\PGF{\C}}^\perp = \Prj{\C}$.
\item The pair $(^{\perp}\GI{\C}, \GI{\C})$ is a complete hereditary cotorsion pair generated by a set such that $^{\perp}\GI{\C}$ is thick and closed under filtered colimits, and $^{\perp}\GI{\C} \cap \GI{\C} = \Inj{\C}$.
\item The pair $(\GF{\C}, {\GF{\C}}^\perp)$ is a complete hereditary cotorsion pair generated by a set such that ${\GF{\C}}^\perp = \Cot{\C} \cap {\PGF{\C}}^\perp$ and $\GF{\C}\cap\PGF{\C}^\perp=\Flat{\C}$, and so $\GF{\C} \cap {\GF{\C}}^\perp = \Flat{\C} \cap \Cot{\C}$.
\end{prt}
\end{lem}

In what follows, we set
$$\E=\{X \in \Mod{\C,A}\ |\ X^{\natural}\in{\PGF{\C}}^\perp\}.$$

The next result gives a projective model structure on $\Mod{\C,A}$, which will then be used to construct a flat model structure.

\begin{prp}\label{prj model}
The pair $(^{\perp}\E, \E)$ is a hereditary projective cotorsion pair in $\Mod{\C,A}$ generated by a set. Moreover, there is a hereditary abelian model structure $(^{\perp}\E, \E, \Mod{\C,A})$ on $\Mod{\C,A}$, where $^{\perp}\E$ is the subcategory of cofibrant objects, $\E$ is the subcategory of trivial objects, and all objects in $\Mod{\C,A}$ are fibrant objects.
\end{prp}

\begin{prf*}
It follows from Lemma \ref{iso}(a) that $\E = \{ G \otimes_\kk A\ |\ G \in \PGF{\C}\}^\perp$, so $(^{\perp}\E, \E)$ is a cotorsion pair in $\Mod{\C,A}$; see \ref{cotpair}. It is clear that $(^{\perp}\E, \E)$ is hereditary as $\E$ is closed under cokernels of monomorphisms. By Lemma \ref{known}(a), the cotorsion pair $(\PGF{\C}, {\PGF{\C}}^\perp)$ is generated by a set $\Se$, that is, $\Se \subseteq \PGF{\C}$ and $\Se^\perp = {\PGF{\C}}^\perp$. One deduces that $\{G\otimes_\kk A\ |\ G\in \Se\}^\perp = \E$ again by Proposition \ref{iso}(a). Consequently, the cotorsion pair $(^{\perp}\E, \E)$ is generated by a set, so $(^{\perp}\E, \E)$ is complete; see e.g. G\"{o}bel and Trlifaj \thmcite[6.11(b)]{RT12}. Since the subcategory ${\PGF{\C}}^\perp$ is thick by Lemma \ref{known}(a), so is $\E$ clearly.

We prove that $\E$ contains the projective objects, which implies that $(^{\perp}\E, \E)$ is a projective cotorsion pair in $\Mod{\C,A}$ by \prpcite[3.7]{Gil16}. For each $X \in \Prj{\C,A}$, one has $\pd_\C X^{\natural}<\infty$ by Proposition \ref{forgive}(a). It follows from Lemma \ref{known}(a) that $\Ext[\C]{i}{G}{P} = 0$ for each $G \in \PGF{\C}$, $P \in \Prj{\C}$ and $i \geqslant 1$, so by dimension shifting one concludes that $X^{\natural}$ is in $\PGF{\C}^\perp$. Thus one has $X\in\E$, as desired.
\end{prf*}

We then give an injective model structure on $\Mod{\C,A}$. Set
\[
\T = \{X \in \Mod{\C,A}\ |\ X^{\natural} \in{^\perp\GI{\C}}\}.
\]

\begin{prp}\label{inj model}
The pair $(\T, \T^\perp)$ is a hereditary cotorsion pair in $\Mod{\C,A}$ such that $\T$ is thick and contains injective objects. If furthermore the subcategory $^\perp\GI{\C}$ is closed under pure quotients, then $(\T, \T^\perp)$ is an injective perfect cotorsion pair. In this case, there is a hereditary abelian model structure $(\Mod{\C,A}, \T, \T^\perp)$ on $\Mod{\C,A}$, where $\T^\perp$ is the subcategory of fibrant objects, $\T$ is the subcategory of trivial objects, and all objects in $\Mod{\C,A}$ are cofibrant objects.
\end{prp}

\begin{prf*}
By Lemma \ref{iso}(b), one has $\T = {^\perp\{\Hom[\kk]{A}{G}\ | \ G \in \GI{\C}\}}$, so $(\T, \T^\perp)$ is a cotorsion pair in $\Mod{\C,A}$; it is hereditary as $\T$ is closed under kernels of epimorphisms. It is clear that the subcategory $\T$ is thick since ${^\perp\GI{\C}}$ is thick by Lemma \ref{known}(b). For each $X \in \Inj{\C,A}$, one has $\id_\C X^{\natural}<\infty$ by Proposition \ref{forgive}(b). It follows from Lemma \ref{known}(b) that $\Ext[\C]{i}{I}{G} = 0$ for each $G \in \GI{\C}$, $I \in \Inj{\C}$ and $i \geqslant 1$, so by dimension shifting one gets that $X^{\natural}$ is in $^\perp\GI{\C}$. Thus one has $X \in \T$, and hence $\T$ contains all injective objects in $\Mod{\C,A}$.

Now we assume that the subcategory $^\perp\GI{\C}$ is closed under pure quotients. Then $^\perp\GI{\C}$ is also closed under pure submodules since $(^{\perp}\GI{\C},\GI{\C})$ is a hereditary cotorsion pair by Lemma \ref{known}(b). The forgetful functor $(-)^{\natural}$ preserves colimits and pure exact sequence by \lemcite[5.7]{HJ21}. It follows that the subcategory $\T$ is closed under pure submodules and pure quotients, and furthermore is closed under filtered colimits since the subcategory $^\perp\GI{\C}$ is closed under filtered colimits by Lemma \ref{known}(b). Thus by \thmcite[A.3]{HJ21}, the cotorsion pair $(\T, \T^\perp)$ is perfect. Moreover, by \prpcite[3.6]{Gil16}, $(\T, \T^\perp)$ is an injective cotorsion pair.
\end{prf*}

In the following we give a flat model structure on $\Mod{\C,A}$. For this purpose, we first prove the next lemma.

\begin{lem}\label{flat cotorsion pair}
The pair $(^{\perp}(\Cot{\C,A}\cap\E), \Cot{\C,A}\cap\E)$ is a complete and hereditary cotorsion pair in $\Mod{\C,A}$ generated by a set, and $^{\perp}(\Cot{\C,A}\cap\E)\cap\E=\Flat{\C,A}$.
\end{lem}

\begin{prf*}
By Proposition \ref{prj model}, there is a set $\Se \subseteq {^\perp\E}$ such that $\Se^\perp = \E$. On the other hand, the cotorsion pair $(\Flat{\C,A}, \Cot{\C,A})$ is generated by a set $\R$ (see \ref{flat-cotorsion}), that is, $\R \subseteq \Flat{\C,A}$ and $\R^\perp = \Cot{\C,A}$. Thus one has
\[
(\R \cup \Se)^\perp = \R^\perp \cap \Se^\perp = \Cot{\C,A} \cap \E,
\]
which implies that $(^{\perp}(\Cot{\C,A} \cap \E), \Cot{\C,A} \cap \E)$ is a complete cotorsion pair in $\Mod{\C,A}$ generated by a set; see \thmcite[6.11(b)]{RT12}. The subcategory $\Cot{\C,A} \cap \E$ is closed under cokernels of monomorphisms as so are $\Cot{\C,A}$ and $\E$ by Proposition \ref{prj model}. Thus the cotorsion pair
\[
(^{\perp}(\Cot{\C,A} \cap \E), \Cot{\C,A} \cap \E)
\]
is hereditary.

We next prove the equality
$^{\perp}(\Cot{\C,A} \cap \E) \cap \E = \Flat{\C,A}$. It is clear that $\Flat{\C,A} = {^{\perp}\Cot{\C,A}} \subseteq {^{\perp}(\Cot{\C,A} \cap \E)}$. On the other hand, if $X \in \Flat{\C,A}$, then $X^{\natural}$ has finite flat dimension in $\Mod{\C}$ by Lemma \ref{forgive}(c), and so it follows from \thmcite[4.11]{SS20} that $X^{\natural}$ is in $\PGF{\C}^\perp$. This implies that $X$ is in $\E$. Thus one has
\[
\Flat{\C,A} \subseteq{^{\perp} (\Cot{\C,A} \cap \E) \cap \E}.
\]
To prove the inclusion of the other direction, we take an arbitrary
$$Y \in{^{\perp}(\Cot{\C,A} \cap \E) \cap\E}.$$
Since $(\Flat{\C,A}, \Cot{\C,A})$ is a complete cotorsion pair (see \ref{flat-cotorsion}), there is an exact sequence
\begin{equation*}\label{eq:ipg1}\tag{\ref{flat cotorsion pair}.1}
0 \to C \to F \to Y \to 0
\end{equation*}
with $C \in \Cot{\C,A}$ and $F \in \Flat{\C,A} \subseteq \E$ as proved before. Thus $C$ is contained in $\E$ as $\E$ is thick by Proposition \ref{prj model}, so $\Ext[\C,A]{1}{Y}{C}=0$. Consequently, the sequence (\ref{eq:ipg1}) splits, and so $Y$ belongs to $\Flat{\C,A}$.
\end{prf*}

Now we are in a position to give the main result in this section.

\begin{thm}\label{flat model structure in general}
There is a hereditary abelian model structure
$$(^{\perp}(\Cot{\C,A}\cap\E), \E, \Cot{\C,A})$$
on $\Mod{\C,A}$, where $^{\perp}(\Cot{\C,A}\cap\E)$ is the subcategory of cofibrant objects, $\E$ is the subcategory of trivial objects, and $\Cot{\C,A}$ is the subcategory of fibrant objects.
\end{thm}
\begin{prf*}
We mention that $(\Flat{\C,A}, \Cot{\C,A})$ is a complete and hereditary cotorsion pair in $\Mod{\C,A}$; see \ref{flat-cotorsion}. Thus by Gillespie \thmcite[1.1]{G15} and Lemma \ref{flat cotorsion pair}, there is a hereditary abelian model structure
$$(^{\perp}(\Cot{\C,A}\cap\E), \W, \Cot{\C,A})$$
on $\Mod{\C,A}$, where the subcategory $\W$ can be described as follows:
\[
\W=\{X\in\Mod{\C,A}\ |\ \exists\ s.e.s\ X\rightarrowtail U\twoheadrightarrow F\ with\ U\in\Cot{\C,A}\cap\E\ and\ F\in\Flat{\C,A}\}.
\]
We then prove that there is an equality $\W=\E$.

Take an object $X\in\W$. Then there is an exact sequence
$$0\to X\to U\to F\to 0$$
in $\Mod{\C,A}$ with $U\in\Cot{\C,A}\cap\E$ and $F\in\Flat{\C,A}$. Thus the sequence
$$0\to X^{\natural}\to U^{\natural}\to F^{\natural}\to 0$$
is exact, where $U^{\natural}\in\PGF{\C}^\perp$ and $F^{\natural}$ has finite flat dimension by Proposition \ref{forgive}(c). So it is easy to check that $F^{\natural}$ is in $\PGF{\C}^\perp$ as $\Flat{\C}\subseteq\PGF{\C}^\perp$ by Lemma \ref{known}(c). The subcategory $\PGF{\C}^\perp$ is thick by Lemma \ref{known}(a), so $X^{\natural}$ is in $\PGF{\C}^\perp$. This implies that $X$ is in $\E$. To prove the inclusion of the other direction, we take an object $Y\in\E$. By the completeness of the cotorsion pair $(\Flat{\C,A}, \Cot{\C,A})$, there is an exact sequence
$$0\to Y \to U' \to F' \to 0$$
in $\Mod{\C,A}$ with $U'\in\Cot{\C,A}$ and $F'\in\Flat{\C,A}$. Consider the exact sequence
$$0\to Y^{\natural} \to (U')^{\natural} \to (F')^{\natural} \to 0$$
in $\Mod{\C}$. The object $(F')^{\natural}$ has finite flat dimension by Proposition \ref{forgive}(c), and so it is in $\PGF{\C}^\perp$ as prove before. On the other hand, $Y^{\natural}$ is in $\PGF{\C}^\perp$ as $Y\in\E$. Thus by Lemma \ref{known}(a), one has $(U')^{\natural}\in\PGF{\C}^\perp$, and so $U'$ is in $\E$. This implies that $Y$ is in $\W$.

Now one gets a hereditary abelian model structure
$$(^{\perp}(\Cot{\C,A}\cap\E), \E, \Cot{\C,A})$$
on $\Mod{\C,A}$, as desired.
\end{prf*}

\begin{setup}\label{setup}
In the rest of this section we let $\Q$ be a small $\kk$-linear category satisfying the following self-dual conditions specified in \cite[Setup 2.5]{HJ21}:
\begin{rqm}
\item Each $\kk$-module $\Q(p,q)$ is finitely generated and projective.
\item $\Q$ satisfies the local boundedness, that is for any object $q$ in $\Q$, there are only finitely many objects in $\Q$ mapping nontrivially into or out of $q$.
\item There exists a Serre functor, that is a $\kk$-linear auto-equivalence $\mathbb{S}: \Q \to \Q$ such that $\Q(p,q)\is\Hom[\kk]{\Q(q,\mathbb{S}(p))}{\kk}.$
\item For any object $q$ in $\Q$, the unit map $\kk\to\Q(q,q)$ (via $x\mapsto x\cdot\id_q$) has a $\kk$-module retraction; whence there is a $\kk$-module decomposition $\Q(q,q)=\kk\cdot\id_q\oplus\tau_q$.
\end{rqm}
\end{setup}

\begin{rmk}\label{show}
If $\kk$ is Gorenstein, then it follows from Dell'Ambrogio, Stevenson and \v{S}\v{t}ov\'{\i}\v{c}ek \thmcite[1.6]{DSS17} that $\Mod{\Q}$ is a locally Gorenstein category in the sense of Enochs, Estrada and Garc\'{\i}a~Rozas \cite{EEG08}, and so it is easy to check that all Gorenstein projective objects in $\Mod{\Q}$ are Gorenstein flat. Thus the subcategory $\PGF{\Q}$ coincides with the subcategory $\GP{\Q}$; see the remark after \exacite[4.10]{SS20}. In this case, Proposition \ref{prj model} recovers \cite[thms. 5.5 and 6.1(a)]{HJ21}.

In Proposition \ref{inj model}, we describe the injective model structure under the assumption that the subcategory $^\perp\GI{\C}$ is closed under pure quotients. We mention that the question whether $^\perp\GI{\C}$ is closed under pure quotients remains open; see \cite[pp. 31]{SS20}. However, for the $\kk$-linear category $\Q$, if $\kk$ is $n$-Gorenstein, then by \thmcite[1.6]{DSS17}, $\Mod{\Q}$ is a locally $n$-Gorenstein category, so $^\perp\GI{\Q}$ is actually the subcategory of objects with finite injective dimension at most $n$. Thus $^\perp\GI{\Q}$ is closed under pure quotients by \cite[prop. 5.3(b) and lem. 5.6]{HJ21}. In this case, Proposition \ref{inj model} recovers \cite[thms. 5.9 and 6.1(b)]{HJ21}.
\end{rmk}

As described in Introduction, the following result gives a flat model structure whose homotopy category is the $\Q$-shaped derived category introduced in \cite{HJ21}.

\begin{cor}\label{flat model structure}
Let $\Q$ be a small $\kk$-linear category satisfying the conditions in Setup \ref{setup}, and assume that $\kk$ is Gorenstein. Then there is a hereditary abelian model structure
\[
(^{\perp}(\Cot{\Q,A} \cap \E), \E, \Cot{\Q,A})
\]
on $\Mod{\Q,A}$, where
$\E = \{X \in \Mod{\Q,A}\ |\ X^{\natural}\ \mathrm{has\ finite\ projective\ dimension\ in}\ \Mod{\Q}\}$.
\end{cor}

\begin{prf*}
The model structure in the statement can be obtained by Theorem \ref{flat model structure in general}. As shown in Remark \ref{show}, $\Mod{\Q}$ is a locally Gorenstein category, so one has
\[
\PGF{\Q}^\perp=\GP{\Q}^\perp=\Prj{\Q}_{<\infty},
\]
where $\Prj{\Q}_{<\infty}$ is the subcategory of objects in $\Mod{\Q}$ with finite projective dimension. Thus one gets that
$$\E = \{X \in \Mod{\Q,A}\ |\ X^{\natural}\ \mathrm{has\ finite\ projective\ dimension\ in}\ \Mod{\Q}\}$$
as desired.
\end{prf*}

\begin{rmk}
If we impose the following extra conditions on $\Q$ and $\kk$:
\begin{enumerate}
\setcounter{enumi}{4}
\item $\kk$ is noetherian and hereditary;
\item $\Q$ satisfies the \emph{strong retraction property}, that is, the $\kk$-submodules $\tau_q$ in Setup \ref{setup}(4) are compatible with composition in $\Q$ as follows:
\begin{center}
$\tau_q \circ \tau_q \subseteq \tau_q$ for all objects $q$, and $\Q(q,p) \circ \Q(p,q) \subseteq\tau_p$ for all $p \neq q$;
\end{center}
\item the pseudo-radical $\tau$ given in \cite[rmk. 7.4 and lem. 7.7]{HJ21} is nilpotent,
\end{enumerate}
then there is an inclusion
\[
^{\perp}(\Cot{\Q,A} \cap \E) \subseteq \setof{X \in \Mod{\Q,A}}{X(q) \in \Flat{A} \text{ for all objects}\ q}.
\]
Indeed, take an arbitrary $X$ in $^{\perp}(\Cot{\Q,A} \cap \E)$, and fix an object $q$ in $\Q$ and an object $N$ in $\Cot{A}$. It follows from \cite[thm. 7.1 and lem. 7.14]{HJ21} that $G_q(N)$ is in $\E$, and hence in $\Cot{\Q,A} \cap \E$ by Corollary \ref{cot}. Thus one gets that
$$\Ext[A]{1}{X(q)}{N} \is \Ext[\C,A]{1}{X}{G_q(N)} = 0$$
by Lemma \ref{n-iso}, and so $X(q)$ belongs to $\Flat{A}$.
%\marginpar{prove the equality as in H-J 2 (need: Gflat, flat)}
\end{rmk}

\section{Characterizations of Gorenstein objects in $\Mod{\Q,A}$}\label{Gorenstein objects}
\noindent
Throughout this section, we assume that $\Q$ is a small $\kk$-linear category satisfying the conditions (1)-(3) in Setup \ref{setup}, and assume that $A$ is an associative $\kk$-algebra (not necessarily has finite projective dimension). If $A$ is Gorenstein then it is known that an object $X$ in $\Mod{\Q, A}$ is Gorenstein projective (resp., Gorenstein injective) if and only if each left $A$-module $X(q)$ is Gorenstein projective (resp., Gorenstein injective) for $q\in\Q$; see \corcite[4.8]{DSS17}. In this section we show that the above result holds for an arbitrary associative $\kk$-algebra $A$. We also prove that an object $X$ in $\Mod{\Q, A}$ is Gorenstein flat (resp., projectively coresolved Gorenstein flat) if and only if each left $A$-module $X(q)$ is Gorenstein flat (resp., projectively coresolved Gorenstein flat) for $q\in\Q$. We prove these results relying on the concept of Frobenius functors introduced by Chen and Ren in \cite{CR22} and a main result of that paper.

\begin{bfhpg}[\bf Frobenius functor]\label{ff}
Let $F: \sf C\to \sf D$ and $G: \sf D \to \sf C$ be two additive functors between abelian categories, and let $\alpha: \sf C \to \sf C$ and $\beta: \sf D \to \sf D$ be two auto-equivalences. Recall from \cite{CR22} that $(F, G)$ is called a \emph{Frobenius pair} (of type $(\alpha, \beta)$) between $\sf C$ and $\sf D$ if both $(F, G)$ and $(G, \beta F\alpha)$ are adjoint pairs, and the functor $F$ is said to be a \emph{Frobenius functor} if it fits into a Frobenius pair $(F, G)$; in this case $G$ is also a Frobenius functor.
\end{bfhpg}

\begin{rmk}\label{3.2}
Let $F: \sf C\to \sf D$ and $G: \sf D \to \sf C$ be two additive functors between Grothendieck categories admitting enough projectives. It is easy to see that if $(F, G)$ is a Frobenius pair then $F$ and $G$ are exact functor preserving projectives, injectives and flats.
\end{rmk}

The next result was proved by Chen and Ren; see \thmcite[3.2]{CR22}.

\begin{lem}\label{CR}
Let $F: \sf C\to \sf D$ be a faithful Frobenius functor between two abelian categories admitting enough projectives. Then an object $X$ in $\sf C$ is Gorenstein projective if and only if $F(X)$ is Gorenstein projective in $\sf D$.
\end{lem}

\begin{ipg}
The product category $\prod_{q\in\Q}\Mod{A}$ is a locally finitely presentable Grothendieck category admitting enough projectives and injectives. Let $(M_{q})_{q\in\Q}$ be an object in $\prod_{q\in\Q}\Mod{A}$. Then it is clear that
$(M_{q})_{q\in\Q}$ is projective (resp., injective, flat) if and only if each left $A$-module $M_{q}$ is projective (resp., injective, flat) for $q\in\Q$. Also, one can check that an object $(M_{q})_{q\in\Q}$ in $\prod_{q\in\Q}\Mod{A}$ is Gorenstein projective\footnote{The definitions of Gorenstein projective, injective and flat objects in the product category $\prod_{q\in\Q}\Mod{A}$ are classical. We mention that for an object $(N_{q})_{q\in\Q}$ in $\prod_{q\in\Q}\rMod{A}$ and an object $(M_{q})_{q\in\Q}$ in $\prod_{q\in\Q}\Mod{A}$, the tensor product is defined as $(N_{q})_{q\in\Q}\otimes(M_{q})_{q\in\Q}=\oplus_{q\in\Q}(N_q\otimes_A M_q)$.} (resp., Gorenstein injective, Gorenstein flat) if and only if each left $A$-module $M_{q}$ is Gorenstein projective (resp., Gorenstein injective, Gorenstein flat) for $q\in\Q$.
\end{ipg}

\begin{lem}\label{Frobenius functor}
There is an adjoint triple $(i^{*},i_{*},i^{!})$ as follows:
\begin{equation*}
  \xymatrix@C=2.5pc{
    \Mod{\Q,A}
    \ar[r]^-{i_{*}}
    &
    \prod_{q\in\Q}\Mod{A}
    \ar@/_1.8pc/[l]_-{i^{*}}
    \ar@/^1.8pc/[l]^-{i^{!}}
  }
  \quad \text{given by} \quad
  {\setlength\arraycolsep{1.5pt}
   \renewcommand{\arraystretch}{1.2}
  \begin{array}{lrl}
  i^{*}(M_{q})_{q\in\Q} &=&\coprod_{q\in\Q} (\Q(q,-) \otimes_\kk M_{q}) \\
  i_{*}(X) &=& (X\mspace{1.5mu}(q))_{q\in\Q} \\
  i^{!}(M_{q})_{q\in\Q} &=&\prod_{q\in\Q} \Hom[\kk]{\Q(-, q)}{M_{q}}\;.
  \end{array}
  }
  \end{equation*}
Moreover, $(i_{*}, i^{!})$ is a Frobenius pair. In particular, $i_{*}$ and $i^{!}$ are Frobenius functors.
\end{lem}
\begin{prf*}
For $X\in\Mod{\Q,A}$ and $(M_{q})_{q\in\Q}\in\prod_{q\in\Q}\Mod{A}$ there are isomorphisms
\begin{align*}
\Hom[\Q,A]{i^{*}(M_{q})_{q\in \Q}}{X}&\cong\prod_{q\in \Q}\Hom[\Q,A]{\Q(q,-)\otimes_{\Bbbk}M_{q}}{X}\\
&\cong\prod_{q\in \Q}\Hom[A]{M_{q}}{X(q)}\\
&\cong\Hom[\prod_{q\in\Q}\Mod{A}]{(M_{q})_{q\in\Q}}{(X(q))_{q\in\Q}}\\
&\cong\Hom[\prod_{q\in\Q}\Mod{A}]{(M_{q})_{q\in\Q}}{i_{*}X},
\end{align*}
where the second isomorphism holds by \ref{adjoint triple}. On the other hand, one has
\begin{align*}
\Hom[\prod_{q\in\Q}\Mod{A}]{i_{*}(X)}{(M_{q})_{q\in \Q}}&\cong\prod_{q\in\Q}\Hom[A]{X(q)}{M_{q}}\\
&\cong\prod_{q\in\Q}\Hom[\Q,A]{X}{\Hom[\Bbbk]{\Q(-,q)}{M_{q}}}\\
&\cong\Hom[\Q,A]{X}{\prod_{q\in\Q}\Hom[\Bbbk]{\Q(-,q)}{M_{q}}}\\
&=\Hom[\Q,A]{X}{i^{!}(M_{q})_{q\in\Q}},
\end{align*}
where the second isomorphism holds again by \ref{adjoint triple}. Thus $(i^{*},i_{*},i^{!})$ is an adjoint triple.

We mention that the Serre functor $\mathbb{S}$ induces an equivalence
$$\mathbb{S}_{*}:\Mod{\Q,A}\to \Mod{\Q,A}$$
with $\mathbb{S}_{*}(X)(q)=X(\mathbb{S}(q))$ and $\mathbb{S}_{*}^{-1}(X)(q)=X(\mathbb{S}^{-1}(q))$ for each $X\in\Mod{\Q,A}$ and $q\in\Q$.
Thus for each $r\in\Q$, one has
\begin{align*}
i^{*}(M_{q})_{q\in\Q}(r)&=\coprod_{q\in\Q}\Q(q,r)\otimes_{\Bbbk}M_{q}\\
&\cong\coprod_{q\in\Q}\Hom[\Bbbk]{\Q(r,\mathbb{S}(q))}{\Bbbk}\otimes_{\Bbbk}M_{q}\\
&\cong\coprod_{q\in\Q}\Hom[\Bbbk]{\Q(r,\mathbb{S}(q))}{M_{q}}\\
&\cong\coprod_{q\in\Q}\Hom[\Bbbk]{\Q(\mathbb{S}^{-1}(r),q)}{M_{q}}\\
&=\prod_{q\in\Q}\Hom[\Bbbk]{\Q(\mathbb{S}^{-1}(r),q)}{M_{q}}\\
&=i^{!}(M_{q})_{q\in\Q}(\mathbb{S}^{-1}(r))\\
&=\mathbb{S}_{*}^{-1}i^{!}(M_{q})_{q\in\Q}(r),
\end{align*}
where the first and third isomorphisms hold by the definition of Serre functors, the second isomorphism holds as $\Q(r,\mathbb{S}(q))$ is a finitely generated and projective $\kk$-module, and the second equation follows from the fact that $\Q$ satisfies the local boundedness.
Thus one has $i^{*}\is\mathbb{S}_{*}^{-1}\circ i^{!}$ by the naturality of the above isomorphisms,  and so $i^{!}\is\mathbb{S}_{*}\circ i^{*}.$
Since $(i^{*},i_{*})$ is an adjoint pair, it is easy to check that
$(i^{!}, i_{*}\circ\mathbb{S}_{*}^{-1})$ is an adjoint pair. This yields that $(i_{*}, i^{!})$ is a Frobenius pair of type $(\mathbb{S}_{*}^{-1}, \id_{\prod_{q\in\Q}\Mod{A}})$.
\end{prf*}

\begin{thm}\label{GP}
An object $X\in \Mod{\Q,A}$ is Gorenstein projective if and only if each left $A$-module $X(q)$ is Gorenstein projective for $q\in\Q$.
\end{thm}
\begin{prf*}
It follows from Lemma \ref{Frobenius functor} that $i_{*}$ is a Frobenius functor, and it is faithful clearly. Thus the desired conclusion holds by Lemma \ref{CR}.
\end{prf*}

Dually, we have the following result for Gorenstein injective objects.

\begin{thm}\label{GI}
An object $X\in \Mod{\Q,A}$ is Gorenstein injective if and only if each left $A$-module $X(q)$ is Gorenstein injective for $q\in\Q$.
\end{thm}

In the following we give a characterization for Gorenstein flat object in $\Mod{\Q,A}$. We mention that by \cite[Corollary 3.9]{HJ21} there is an adjoint triple $(F'_{q},E'_{q},G'_{q})$ for each $q\in\Q$ as follows (see \ref{adjoint triple}):
  \begin{equation*}
  \xymatrix@C=4pc{
    \rMod{\Q,A}
    \ar[r]^-{E'_{q}}
    &
    \rMod{A}
    \ar@/_1.8pc/[l]_-{F'_{q}}
    \ar@/^1.8pc/[l]^-{G'_{q}}
  }
  \qquad \text{given by} \qquad
  {\setlength\arraycolsep{1.5pt}
   \renewcommand{\arraystretch}{1.2}
  \begin{array}{rcl}
  F'_q(M) &=& \Q(-,q) \otimes_\kk M \\
  E'_q(X) &=& X\mspace{1.5mu}(q) \\
  G'_{q}(M) &=& \Hom[\kk]{\Q(q, -)}{M}\;.
  \end{array}
  }
  \end{equation*}

\begin{lem}\label{technical}
For each $M\in \Mod{A}$, $N\in \rMod{A}$ and $E\in \Mod{\kk}$, there are isomorphisms
$$G'_{q}(\Hom[\Bbbk]{M}{E})\cong\Hom[\Bbbk]{G_{\mathbb{S}(q)}(M)}{E}$$
and
$$G_{\mathbb{S}(q)}(\Hom[\Bbbk]{N}{E})\cong\Hom[\Bbbk]{G'_{q}(N)}{E}.$$
\end{lem}
\begin{prf*}
We only prove the first isomorphism; the second one is proved similarly.

For each $p\in\Q$, one has
\begin{align*}
G'_{q}(\Hom[\Bbbk]{M}{E})(p)&=\Hom[\Bbbk]{\Q(q,p)}{\Hom[\Bbbk]{M}{E}}\\
&\cong\Hom[\Bbbk]{\Q(q,p)\otimes_{\Bbbk}M}{E}\\
&\cong\Hom[\Bbbk]{\Hom[\Bbbk]{\Q(p,\mathbb{S}(q))}{\Bbbk}\otimes_{\Bbbk}M}{E}\\
&\cong\Hom[\Bbbk]{\Hom[\Bbbk]{\Q(p,\mathbb{S}(q))}{M}}{E}\\
&=\Hom[\Bbbk]{G_{\mathbb{S}(q)}(M)}{E}(p),
\end{align*}
where the second isomorphism holds as $\mathbb{S}$ is a Serre functor, and the third one holds since the $\kk$-module $\Q(p,\mathbb{S}(q))$ is finitely generated and projective. It follows that $G'_{q}(\Hom[\Bbbk]{M}{E})\cong\Hom[\Bbbk]{G_{\mathbb{S}(q)}(M)}{E}$ by the naturality of the above isomorphisms.
\end{prf*}

\begin{lem}\label{key formula}
For each $N\in \rMod{A}$ and $X\in\Mod{\Q,A}$, there is an isomorphism
$$G'_{q}(N)\otimes_{Q,A}X\is N\otimes_{A}X(\mathbb{S}(q)).$$
\end{lem}

\begin{prf*}
For each $E\in \Mod{\Bbbk}$, there are isomorphisms
\begin{align*}
\Hom[\Bbbk]{G'_{q}(N)\otimes_{\Q,A}X}{E}&\cong\Hom[\Q,A]{X}{\Hom[\Bbbk]{G'_{q}(N)}{E}}\\
&\cong\Hom[\Q,A]{X}{G_{\mathbb{S}(q)}(\Hom[\Bbbk]{N}{E})}\\
&\cong\Hom[A]{X(\mathbb{S}(q))}{\Hom[\Bbbk]{N}{E}}\\
&\cong\Hom[\Bbbk]{N\otimes_{A}X(\mathbb{S}(q))}{E},
\end{align*}
where the first isomorphism follows from (\ref{1.6.1}), the second one holds by Lemma \ref{technical}, and the third one holds as $(E_{\mathbb{S}(q)},G_{\mathbb{S}(q)})$ is an adjoint pair; see \ref{adjoint triple}. Thus there is an isomorphism $G'_{q}(N)\otimes_{Q,A}X\is N\otimes_{A}X(\mathbb{S}(q))$.
\end{prf*}

\begin{lem}\label{lem3.9}
Let $X$ be an object in $\Mod{\Q,A}$ such that each left $A$-module $X(q)$ is Gorenstein flat for $q\in\Q$. Then $\Tor[i]{\Q,A}{E}{X}=0$ for each injective object $E$ in $\rMod{\Q,A}$ and all $i\geq 1$.
\end{lem}
\begin{prf*}
Fix an injective object $E$ in $\rMod{\Q,A}$. Let $I$ be a faithful injective right $A$-module.
Then for all $q\in\Q$ the objects $G'_{q}(I)$ are injective and cogenerate $\lMod{\Q,A}$; see \cite[Proposition 3.12]{HJ21}.
Thus $E$ is a direct summand of an object of the form $\prod_{q\in\Q}\ G'_{q}(I).$
We mention that
$\prod_{q\in\Q}\ G'_{q}(I) \cong \coprod_{q\in\Q} G'_{q}(I)$ by \cite[Proposition 3.7]{HJ23}.
Since the functor $\Tor[i]{\Q,A}{-}{X}$ converts coproducts to coproducts, without loss of generality we may assume that $E$ has the form $G'_{q}(I)$. Let $\mathbf{P}$ be a projective resolution of $X$. Then for each $i\geq 1$ one has
\begin{align*}
\Tor[i]{\Q,A}{G'_{q}(I)}{X}&\cong \mathrm{H}_{i}(G'_{q}(I)\otimes_{\Q,A}\mathbf{P})\\
&\cong \mathrm{H}_{i}(I\otimes_{A}\mathbf{P}(\mathbb{S}(q)))\\
&\cong\Tor[i]{A}{I}{X(\mathbb{S}(q))}\\
&=0,
\end{align*}
where the second isomorphism holds by Lemma \ref{key formula} and the equality holds as $X(\mathbb{S}(q))$ is Gorenstein flat. This completes the proof.
\end{prf*}

\begin{thm}\label{GF}
An object $X\in \Mod{\Q,A}$ is Gorenstein flat if and only if each left $A$-module $X(q)$ is Gorenstein flat for $q\in\Q.$
\end{thm}
\begin{prf*}
For the ``only if" part, let $X$ be a Gorenstein flat object in $\Mod{\Q,A}$. Then there is an exact sequence
$$\mathbf{F}= \cdots\to F_{1}\to F_{0}\to F_{-1}\to F_{-2}\to\cdots$$
of flat objects in $\Mod{\Q,A}$ such that $X\cong \Coker(F_{1}\to F_{0})$ and the sequence remains exact after applying the functor $E\otimes_{\Q,A}-$ for each injective object $E$ in $\rMod{\Q,A}$.
Let $I$ be an injective right $A$-module.
Then it follows from Lemma \ref{key formula} that there is an isomorphism
$I\otimes_{A}\mathbf{F}(q)\cong G'_{\mathbb{S}^{-1}(q)}(I)\otimes_{\Q,A}\mathbf{F}$ for each $q\in\Q$.
We mention that $G'_{\mathbb{S}^{-1}(q)}(I)$ is an injective object in $\rMod{\Q,A}$.
Then the $\kk$-complex $G'_{\mathbb{S}^{-1}(q)}(I)\otimes_{\Q,A}\mathbf{F}$ is exact, and so $I\otimes_{A}\mathbf{F}(q)$ is exact.
On the other hand, $\mathbf{F}(q)$ is an exact sequence of flat left $A$-modules by Proposition \ref{degreewise}, which yields that $X(q)$ is Gorenstein flat for each $q\in\Q$.

For the ``if" part, let $X$ be an object in $\Mod{\Q,A}$ such that each left $A$-module $X(q)$ is Gorenstein flat for $q\in\Q$. Then $i_{*}(X)=(X(q))_{q\in\Q}$ is a Gorenstein flat object in the product category $\prod_{q\in\Q}\lMod{A}$.
Thus there is an exact sequence
$$0\to i_{*}(X)\xrightarrow{f'}F\to Y\to 0$$
in $\prod_{q\in\Q}\lMod{A}$ with $F$ flat and $Y$ Gorenstein flat.
Denote by $\eta$ and $\varepsilon$ the unit and counit of the adjoint pair $(i_{*},i^{!})$, respectively.
Set $f=i^{!}(f')\circ \eta_{X}:X\to i^{!}(F)$.
Since $i_{*}$ is faithful and $i^{!}$ is exact, it follows that $f$ is a monomorphism.
Consequently, one gets an exact sequence
$$0\to X\xrightarrow{f}i^{!}(F)\to X'\to 0$$
in $\lMod{\Q,A}$.
Since $f'=\varepsilon_{F}\circ i_{*}(f)$, there is a pull-back diagram
$$\xymatrix{0\ar[r]&i_{*}(X)\ar[r]^{i_{*}(f)}\ar@{=}[d]_{}&i_{*}i^{!}(F)\ar[r]^{}\ar[d]_{\varepsilon_{F}}&i_{*}(X')\ar[r]\ar[d]&0\\
    0\ar[r]&i_{*}(X)\ar[r]^{f'}&F\ar[r]^{}&Y \ar[r]&0}$$
    which yields an exact sequence
    $$0\to i_{*}i^{!}(F)\to i_{*}(X')\oplus F\to Y\to 0.$$
We mention that $i_{*}$ and $i^{!}$ preserves flats as $(i_{*}, i^{!})$ is a Frobenius pair by Lemma \ref{Frobenius functor}. Then it follows from \cite[Corollary 4.12]{SS20} that $i_{*}(X')$ is a Gorenstein flat object in the product category $\prod_{q\in \Q}\lMod{A}$. Thus for each injective object $E$ in $\rMod{\Q,A}$ one has $\Tor[i]{\Q,A}{E}{X'}=0$ for all $i\geq 1$ by Lemma \ref{lem3.9}.
Set $F_{0}=i^{!}(F)$ and repeat the above argument for $X'$.
One obtains inductively an exact sequence
$$0\to X\to F_{0}\xrightarrow{d_{0}} F_{-1}\xrightarrow{d_{-1}} F_{-2}\to \cdots$$
in $\lMod{\Q,A}$ with each $F_{i}$ flat such that it remains exact after applying the functor $E\otimes_{\Q,A}-$ for every injective object $E$ in $\rMod{\Q,A}$.
Thus it follows from Lemma \ref{lem3.9} that $X$ is an Gorenstein flat object in $\lMod{\Q,A}$.
\end{prf*}

We end this section with the following result for projectively coresolved Gorenstein flat objects, which can be proved similarly as in Theorem \ref{GF}.

\begin{thm}\label{PGF}
An object $X\in \Mod{\Q,A}$ is projectively coresolved Gorenstein flat if and only if each left $A$-module $X(q)$ is projectively coresolved Gorenstein flat for $q\in\Q.$
\end{thm}

\section*{Acknowledgment}

\noindent
We thank Henrik Holm and Peter J{\o}rgensen for valuable comments on a draft of this paper. Thanks are also due to the anonymous referee for helpful suggestions that improved the exposition.

\bibliographystyle{amsplain-nodash}
%\bibliographystyle{amsplain}
%\bibliography{../+references}

\def\cprime{$'$}
  \providecommand{\arxiv}[2][AC]{\mbox{\href{http://arxiv.org/abs/#2}{\sf
  arXiv:#2 [math.#1]}}}
  \providecommand{\oldarxiv}[2][AC]{\mbox{\href{http://arxiv.org/abs/math/#2}{\sf
  arXiv:math/#2
  [math.#1]}}}\providecommand{\MR}[1]{\mbox{\href{http://www.ams.org/mathscinet-getitem?mr=#1}{#1}}}
  \renewcommand{\MR}[1]{\mbox{\href{http://www.ams.org/mathscinet-getitem?mr=#1}{#1}}}
\providecommand{\bysame}{\leavevmode\hbox to3em{\hrulefill}\thinspace}
\providecommand{\MR}{\relax\ifhmode\unskip\space\fi MR }
\providecommand{\MRhref}[2]{%
  \href{http://www.ams.org/mathscinet-getitem?mr=#1}{#2}
}
\providecommand{\href}[2]{#2}

\end{document}